\documentclass[12pt]{article}

\usepackage{dcolumn}
\usepackage{amsmath}
\usepackage{amssymb}

\makeatletter
\def\btt#1{\texttt{\@backslashchar#1}}%
\DeclareRobustCommand\bblash{\btt{\@backslashchar}}%
\makeatother
\def\b#1{{\mathbb #1}}
\def\c#1{{\cal #1}}

\def\1{{\bf 1}}

\def\nn{\nonumber \\}

\def\sq{\mbox{\rlap{$\sqcap$}$\sqcup$}}

\def\R{{\cal R}\,}

\newcommand{\tl}{\,\triangleleft}

\def\cocross{{>\!\!\!\triangleleft}}
\def\g{\mbox{\bf g\,}}
\def\uqg{\mbox{$U_q{\/\mbox{\bf g}}$ }}
\def\uqs{\mbox{$U_q{\/\mbox{so}(N)}$}}

\newcommand{\be}{\begin{equation}}
\newcommand{\ee}{\end{equation}}
\newcommand{\bea}{\begin{eqnarray}}
\newcommand{\eea}{\end{eqnarray}}
\newcommand{\ba}{\begin{array}}
\newcommand{\ea}{\end{array}}
\newtheorem{prop}{Proposition}

\newtheorem{theorem}{Theorem}
\newtheorem{corollary}{Corollary}
\newenvironment{proof}[1]{\vspace{5pt}\noindent{\bf Proof #1}\hspace{6pt}}%
{\hfill\sq}
\newcommand{\bp}{\begin{proof}}
\newcommand{\ep}{\end{proof}\par\vspace{10pt}\noindent}
\begin{document}

\title{Decoupling Braided Tensor Factors
\footnote{Talk given at the 23-rd International 
Conference on Group Theory Methods in Physics, Dubna (Russia), August 2000}}
\author{Gaetano Fiore,$\strut^{1,2}$ \, Harold Steinacker$\strut^{3}$,
        Julius Wess$\strut^{3,4}$ \\\\
        \and
        $\strut^1$Dip. di Matematica e Applicazioni, Fac.  di Ingegneria\\ 
        Universit\`a di Napoli, V. Claudio 21, 80125 Napoli
        \and
        $\strut^2$I.N.F.N., Sezione di Napoli,\\
        Complesso MSA, V. Cintia, 80126 Napoli
        \and
        $\strut^3$Sektion Physik, Ludwig-Maximilian Universit\"at,\\
        Theresienstra\ss e 37, D-80333 M\"unchen
        \and
        $\strut^4$Max-Planck-Institut f\"ur Physik\\
        F\"ohringer Ring 6, D-80805 M\"unchen
        }

\date{}
\maketitle

\begin{abstract}
We briefly report on our result \cite{FioSteWes00}
that the braided tensor product algebra of two module
algebras $\c{A}_1, \c{A}_2$ of a quasitriangular Hopf algebra $H$ 
is equal to the ordinary tensor product algebra of $\c{A}_1$ with
a subalgebra 
isomorphic to $\c{A}_2$ and {\it commuting} with $\c{A}_1$, 
provided there exists a realization of $H$ within $\c{A}_1$. 
As applications of the theorem we consider the braided tensor product 
algebras of two or more quantum group covariant quantum spaces or deformed 
Heisenberg algebras.
\end{abstract}

\section{Introduction and main theorem\protect}
\label{sec:level1}

As is well known, given two associative unital algebras 
$\c{A}_1, \c{A}_2$ (over the field $\b{C}$, say), there is an
obvious way to build a new algebra $\c{A}$  which is 
as a vector space the tensor product 
$\c{A}=\c{A}_1\otimes_{\b{C}}\c{A}_2$
of the two vector spaces (over the same field) and has
a product law such 
that $\c{A}_1\otimes\1$  and $\1\otimes\c{A}_2$ are subalgebras
isomorphic to  $\c{A}_1$ and $\c{A}_2$ respectively:
one just completes the product law by
postulating the trivial commutation relations
\begin{equation}
(\1\otimes a_2)(a_1\otimes \1)=(a_1\otimes \1)(\1\otimes a_2) 
\label{trivial4}
\end{equation}
for any $a_1\in\c{A}_1$, $a_2\in\c{A}_2$. The resulting algebra
is the ordinary tensor product algebra. 
With a standard abuse of notation we shall denote in the sequel
$a_1\otimes a_2$ by $a_1a_2$ for any $a_1\in\c{A}_1$, 
$a_2\in\c{A}_2$; consequently (\ref{trivial4}) becomes 
\begin{equation}
a_2a_1=a_1a_2.                                       \label{trivial}
\end{equation}

If $\c{A}_1, \c{A}_2$ are module algebras of a Lie algebra
$\g$, and we require $\c{A}$ to be too, then (\ref{trivial})
has no alternative, because any 
$g\in\g$ acts as a derivation on the (algebra as well as tensor)
product of any two elements, or,  in Hopf algebra language,
because the coproduct  $\Delta(g)=g_{(1)}\otimes g_{(2)}$ 
(at the rhs we have used Sweedler notation)
of the Hopf algebra $H\equiv U\g$ is cocommutative.
In this paper we shall work with 
right-module algebras (instead of left ones), and denote by 
$\tl:(a_i,g)\in\c{A}_i\times H\rightarrow a_i\tl g\in\c{A}_i$ 
the right action; the reason is that they are equivalent 
to left comodule algebras, which 
are used in much of the literature. In Ref. \cite{FioSteWes00}
we give also the corresponding
formulae for the left module algebras.
We recall that a right action 
$\tl:(a,g)\in\c{A}\times H\rightarrow a\tl g\in\c{A}$ 
by definition fulfills
\begin{eqnarray}
&&a\tl(gg') = (a\tl g) \tl   g',                     \label{modalg1}\\
&&(aa')\tl g = (a\tl g_{(1)})\, (a'\tl g_{(2)}).       \label{modalg2}
\end{eqnarray}

If we take as Hopf algebra $H$ a quasitriangular noncocommutative one
like the quantum group $\uqg$, as
$\c{A}_i$ some $H$-module algebras, and
we require $\c{A}$ to be a $H$-module algebra too, then
(\ref{trivial}) has to be replaced by  one of the formulae
\begin{eqnarray}
&&a_2a_1= (a_1\tl \R^{(1)})\, (a_2\tl \R^{(2)}),     \label{absbraiding}\\
&&a_2a_1= (a_1\tl \R^{-1(2)}) \, (a_2\tl \R^{-1(1)}). \label{absbraiding'}
\end{eqnarray}
This yields instead of $\c{A}$ two different {\it braided} tensor product 
algebras \cite{JoyStr86,Maj1}, which we shall call
$\c{A}^+=\c{A}_1\underline{\otimes}^+\c{A}_2$ and
$\c{A}^-=\c{A}_1\underline{\otimes}^-\c{A}_2$ respectively.
Here $\R\equiv \R^{(1)}\otimes \R^{(2)}\in H^+\otimes H^-$ 
denotes the so-called universal $R$-matrix of $H\equiv$ 
\cite{Dri86}, $\R^{-1}$ its inverse,
and $H^{\pm}$ denote the Hopf positive and negative 
Borel subalgebras of $H$. If in particular $H$ is
triangular, then $\R^{-1}=\R_{21}$, $\c{A}^+=\c{A}^-$,
and one has just one braided tensor product algebra.
In any case, both $\c{A}^+$ and $\c{A}^-$ go to the ordinary 
tensor product algebra
$\c{A}$ in the limit $q\to 1$, because in this limit
$\R\to \1\otimes\1$.

The braided tensor product 
is a particular example of a more general notion, that of a
{\it crossed (or twisted) tensor product} \cite{BorMar00}
of two unital associative algebras.

In view of (\ref{absbraiding}) or (\ref{absbraiding'}) 
studying representations of $\c{A}^{\pm}$ is a more difficult
task than just studying the representations of $\c{A}_1,\c{A}_2$
and taking their tensor products. The degrees of freedom
of $\c{A}_1,\c{A}_2$ are so to say ``coupled''. One might ask 
whether one can ``decouple'' them by a transformation of generators.
As shown in Ref. \cite{FioSteWes00}, the answer is positive if
there respectively exists an algebra homomorphism $\varphi_1^+$ or
an algebra homomorphism $\varphi_1^-$
\begin{equation}
\varphi_1^{\pm}: \c{A}_1\cocross H^{\pm} \rightarrow \c{A}_1
                                            \label{Hom}
\end{equation}
acting as the identity on $\c{A}_1$, namely for any $a_1\in \c{A}_1$
\begin{equation}
\varphi^{\pm}_1(a_1)=a_1  .                          \label{ident0}
\end{equation}
(Here $\c{A}_1\cocross H^{\pm}$ denotes the cross product between 
$\c{A}_1$ and $H^{\pm}$). In other words, this amounts
to assuming that $\varphi_1^+(H^+)$ [resp. $\varphi_1^-(H^-)$] provides a
{\it realization} of $H^+$ (resp. $H^-$) within $\c{A}_1$.
In this report we summarize the main results of Ref. 
\cite{FioSteWes00}. The basic one is

\begin{theorem} \cite{FioSteWes00}.
Let $\{H,\R\}$ be a quasitriangular Hopf algebra and 
$H^+,H^-$ be Hopf subalgebras of $H$ such that $\R\in H^+\otimes H^-$.
Let $\c{A}_1, \c{A}_2$ be respectively a $H^+$- and
a $H^-$-module algebra, so that we can define $\c{A}^+$ 
as in (\ref{absbraiding}), and $\varphi_1^+$ be a homomorphism
of the type (\ref{Hom}), (\ref{ident0}), so that we can define
the map $\chi^+:\c{A}_2\rightarrow \c{A}^+$ by
\begin{equation}
\chi^+(a_2):= \varphi_1^+(\R^{(1)})\, (a_2\tl \R^{(2)}). 
\label{Def+}
\end{equation}
Alternatively, let $\c{A}_1, \c{A}_2$ be respectively a $H^-$- and
a $H^+$-module algebra, so that we can define $\c{A}^-$ 
as in (\ref{absbraiding'}), and $\varphi_1^-$ be a homomorphism
of the type (\ref{Hom}), (\ref{ident0}), so that we can define
the map $\chi^-:\c{A}_2\rightarrow \c{A}^-$ by
\begin{equation}
\chi^-(a_2):= \varphi_1^-(\R^{-1(2)})\, (a_2\tl \R^{-1(1)}). 
\label{Def-}
\end{equation}
In either case $\chi^{\pm}$ are then injective algebra homomorphisms and 
\begin{equation}
[\chi^{\pm}(a_2),\c{A}_1]=0,                       \label{commu}
\end{equation}
namely the subalgebras 
$\tilde\c{A}_2^{\pm}:=\chi^{\pm}(\c{A}_2)\approx\c{A}_2$ commute
with $\c{A}_1$. Moreover $\c{A}^{\pm}=\c{A}_1\otimes\tilde\c{A}_2^{\pm}$.
\label{maintheo}
\end{theorem}

The last equality means that $\c{A}^{\pm}$ are respectively
equal to the ordinary tensor product algebra of $\c{A}_1$ with the
subalgebras $\tilde\c{A}_2^{\pm}\subset\c{A}^{\pm}$, which are isomorphic 
to $\c{A}_2$! $\chi^+,\chi^-$ will be called "unbraiding" maps.

We recall the content of 
the hypotheses stated in the theorem.
The algebra $\c{A}_1\cocross H^{\pm}$ as a vector space is the tensor
product $\c{A}_1\otimes_{\b{C}}H^{\pm}$, as an algebra it has subalgebras
$\c{A}_1\otimes\1$, $\1\otimes H$ and has
cross commutation relations
\begin{equation}
a_1 g=g_{(1)}\, (a_1\tl g_{(2)}),                          \label{crossprod}
\end{equation}
for any $a_1\in\c{A}_1$ and $g\in H^{\pm}$.
$\varphi_1^{\pm}$ being an algebra homomorphism means that for any 
$\xi,\xi'\in\c{A}_1\cocross H^{\pm}$
$\varphi^{\pm}_1(\xi\xi')=\varphi^{\pm}_1(\xi)\,\varphi^{\pm}_1(\xi')$.      
Applying $\varphi^{\pm}_1$ to both sides of (\ref{crossprod})
we find
$a\varphi^{\pm}(g)=\varphi^{\pm}(g_{(1)}) (a\tl g_{(2)})$.

Of course,
we can use the above theorem  iteratively to completely
unbraid the  braided tensor product algebra of an arbitrary number $M$ of 
copies of $\c{A}_1$.  We end up with

\begin{corollary} If $\c{A}_1$ is a (right-) module algebra of the Hopf 
algebra $H$ and there exists an algebra homomorphism $\varphi_1^+$ 
of the type (\ref{Hom}), (\ref{ident0}), then there is an algebra
isomorphism.
\begin{equation}
\underbrace{\c{A}_1\underline{\otimes}^+...\underline{\otimes}^+\c{A}_1}
\limits_{M\mbox{ times}}\approx\underbrace{\c{A}_1\otimes...\otimes\c{A}_1}
\limits_{M\mbox{ times}}.
\end{equation}
An analogous claim holds for the second  braided tensor product
if there exists a map $\varphi_1^-$.
\end{corollary}

\section{The unbraiding under the $*$-structures}
\label{star}

$\c{A}^+$ (as well as $\c{A}^-$) is a $*$-algebra 
if $H$ is a Hopf $*$-algebra, $\c{A}_1$, $\c{A}_2$ are $H$-module
$*$-algebras (we shall use the same symbol $*$ for the $*$-structure
on all algebras $H,\c{A}_1$, etc.), and 
\begin{equation}
\R^* =\R^{-1}                                       \label{real1}
\end{equation}
(here $\R^* $ means $\R^{(1)}{}^*\otimes \R^{(2)}{}^*$).
In the quantum group case (\ref{real1}) requires $|q|=1$.
Under the same assumptions also $\c{A}_1\cocross H$ is
a $*$-algebra. If $\varphi_1^{\pm}$ exist
setting $\varphi_1'{}^{\pm}:=*\circ\varphi_1^{\pm}\circ *$
we realize that also $\varphi_1'{}^{\pm}$ are algebra homomorphisms
of the type (\ref{Hom}), (\ref{ident0}).
If such homomorphisms are uniquely determined, we conclude
that $\varphi_1{}^{\pm}$ are $*$-homomorphisms. More generally,
one may be able to choose $\varphi_1{}^{\pm}$ 
as $*$-homomorphisms.  How do the corresponding $\chi^{\pm}$
behave under $*$?

\begin{prop}\cite{FioSteWes00}. Assume that the conditions of Theorem
\ref{maintheo} for defining $\chi^+$ (resp. $\chi^-$) are
fulfilled. If $\R^*=\R^{-1}$ and $\varphi_1^+$ (resp. $\varphi_1^-$)
is a $*$-homomorphism then $\chi^+$ (resp. $\chi^-$) is, too.
Consequently, $\c{A}_1$, $\tilde\c{A}_2^{\pm}$ are closed under $*$.
\label{propstar}
\end{prop}

\section{Applications}
\label{applications}

In this section we illustrate the application of
Theorem \ref{maintheo} and Corollary 1 
to some algebras $H,\c{A}_i$ for which 
homomorphisms $\varphi_1^{\pm}$ are known. $H$
will be the quantum group $U_qsl(N)$ or $\uqs$, and 
$\c{A}_1$ is the  $U_qsl(N)$- or \uqs-covariant Heisenberg algebra
(Section \ref{heisenberg}), the $\uqs$-covariant
quantum space/sphere (Section \ref{qspaces}). In Ref. \cite{FioSteWes00}
we have treated also the $U_qso(3)$-covariant $q$-fuzzy sphere. 
As generators of H it will be convenient in either case to use 
the Faddeev-Reshetikhin-Takhtadjan (FRT) 
generators \cite{FadResTak89} 
$\c{L}^+{}_l^a\in H^+$ and $\c{L}^-{}_l^a\in H^-$. They are related to $\R$ by
\begin{equation}
\c{L}^+{}_l^a:=\R^{(1)} \rho_l^a(\R^{(2)})\qquad\qquad
\c{L}^-{}_l^a:=\rho_l^a(\R^{-1}{}^{(1)})\R^{-1}{}^{(2)}, \label{frt}
\end{equation}
where $\rho_l^a(g)$ denote the matrix elements of $g\in\uqg$
in the fundamental $N$-dimensional representation $\rho$
of $\uqg$.
In fact they provide, together with the square roots of the
elements $\c{L}^{\pm}{}^i_i$,  a (overcomplete)
set of generators of \uqg.

\subsection{Unbraiding `chains' of braided Heisenberg algebras}
\label{heisenberg}

In this subsection we consider the braided tensor product of
$M\ge 2$ copies of the
$\uqg$-covariant deformed Heisenberg algebras 
$\c{D}_{\epsilon,\g}$, $\g=sl(N),so(N)$. Such algebras 
have been introduced in Ref. 
\cite{PusWor89,WesZum90,CarSchWat91}. They are
unital associative algebras generated by $x^i,\partial_j$
fulfilling the relations
\begin{equation}
\c{P}_a{}^{ij}_{hk}x^hx^k=0, \qquad
\c{P}_a{}^{ij}_{hk}\partial_j\partial_i=0,
\qquad\partial_i x^j=\delta^i_j+(q\gamma\hat R)^{\epsilon}
{}_{ih}^{jk}x^h\partial_k,
\label{gringo}
\end{equation}
where $\gamma=q^{\frac 1N},1$ respectively for $\g=sl(N),so(\!N\!)$,
and the exponent $\epsilon$ can take either value $\epsilon=1,-1$.
$\hat R$ denotes the braid matrix of $\uqg$ [given in formulae
(\ref{defRslN})], and the matrix
$\c{P}_a$ is the deformed
antisymmetric projector appearing in the decomposition 
(\ref{projectorR}) of the latter. 
The coordinates $x^i$ transform according to the fundamental
$N$-dimensional representation $\rho$ of $\uqg$, 
whereas the `partial derivatives' transform according the
contragredient representation,
\begin{equation}
x^i\tl g=\rho^i_j(g)x^j,   \qquad\qquad
\partial_i\tl g=
\partial_h\rho^h_i(S^{-1}g).                  \label{fund}
\end{equation}
In our conventions
the indices will take the values  $i=1,...,N$ if $\g=sl(N)$,
whereas if $\g=so(N)$ they will take the values
$i=-n,\ldots,-1,0,1,\ldots n$ for $N$ odd,
and $i=-n,\ldots,-1, 1,\ldots n$ for $N$ even; 
here $n:=\left[\frac N 2\right]$ 
denotes the rank of $so(N)$. We shall
enumerate the different copies of $\c{D}_{\epsilon,\g}$ by attaching
to them an additional Greek index, e.g. $\alpha=1,2,...,M$. 
The prescription (\ref{absbraiding'}) 
gives the following ``cross'' commutation relations 
between their respective generators ($\alpha<\beta$).
\begin{equation}
\ba{ll}
x^{\alpha,i}x^{\beta,j}=\hat R^{ij}_{hk}x^{\beta,h}x^{\alpha,k},
&\quad\quad\quad\partial_{\alpha,i}\partial_{\beta,j}=
\hat R_{ji}^{kh}\partial_{\beta,h}\partial_{\alpha,k},\\
\partial_{\alpha,i}x^{\beta,j}=
\hat R^{-1}{}^{jh}_{ik}x^{\beta,k}\partial_{\alpha,h},&\quad\quad\quad
\partial_{\beta,i}x^{\alpha,j}=
\hat R^{jh}_{ik}x^{\alpha,k}\partial_{\beta,h}.
\ea
\label{cross}
\end{equation}

Algebra homomorphisms
$\varphi_1:\c{A}_1\cocross H\rightarrow \c{A}_1$, for 
$H=\uqg$  and $\c{A}_1$ equal to
(a suitable completion of) $\c{D}_{\epsilon,\g}$
have been constructed in Ref. \cite{Fiocmp95,ChuZum95}.
This is the $q$-analog of the well-known fact that the elements of  
$\g$ can be realized as ``vector fields''
(first order differential operators) on the corresponding
$\g$-covariant (undeformed) space, e.g.
$\varphi_1(E^i_j)=x^i\partial_j-\frac 1N \delta^i_j$
in the $\g=sl(N)$ case.
The maps $\varphi_1^{\pm}$ needed to apply Theorem \ref{maintheo}
are simply the restrictions to $\c{A}_1\cocross H^{\pm}$
of $\varphi_1$ of Ref. \cite{Fiocmp95,ChuZum95}.

The unbraiding procedure is recursive. We just describe the first step,
which consists of using the
homomorphism $\varphi_1^{\pm}$ 
to unbraid the first copy from the others.
According to the main theorem, if we set 
\begin{eqnarray}
y^{1,i} &\equiv& x^{1,i}        \qquad\qquad\qquad\qquad
\partial_{y,1,a} \equiv \partial_{1,a}         \label{def2.1}\\
y^{\alpha,i} &\equiv& \chi^-(x^{\alpha,i})=
\varphi_1(\R^{-1(2)})\rho^i_j(\R^{-1(1)}) x^{\alpha,j} =
\varphi_1(\c{L}^-{}^i_j)x^{\alpha,j},  \label{def2.3}\\
\partial_{y,\alpha,a} &\equiv& \chi^-(\partial_{\alpha,a})=
\varphi_1(S\R^{-1(2)})\rho^d_a(\R^{-1(1)}) \partial_{\alpha,d}=
\varphi_1(S\c{L}^-{}^d_a) \partial_{\alpha,d},     \label{def2.4}
\end{eqnarray}
with $\alpha>1$. By Theorem \ref{maintheo} 
$y^{1,i}\equiv x^{1,i}$ and $\partial_{y,1,i}\equiv \partial_{1,i}$
will commute with 
$y^{2,i},...,y^{M,i}$ and  $\partial_{y,2,i},...,\partial_{y,M,i}$.
As we see, the FRT generators are special because they
appear in the redefinitions (\ref{def2.3}-\ref{def2.4}).
The explicit expression of $\varphi_1(\c{L}^-{}^i_j)$
in terms of $x^{1,i},\partial_{1,a}$ for $U_qsl(2),U_qso(3)$ 
has been given in Ref. \cite{FioSteWes00}. For different values
of $N$ it can be found from the results of Ref. 
\cite{Fiocmp95,ChuZum95} by passing from the generators
adopted there to the FRT generators.

By completely analogous arguments one
determines the alternative unbraiding procedure 
for the braided tensor product stemming
from prescription (\ref{absbraiding}). 

$\c{A}_1\cocross H$ is a $*$-algebra and
the map $\varphi_1$ is a $*$-homomorphism both for
$q$ real and $|q|=1$. But $\varphi_1^{\pm}$ are $*$-homomorphisms
only for $|q|=1$. In the latter case
 the $*$-structure of $\c{A}_1$ is
\begin{equation}
(x^i)^*=x^i, \qquad (\partial_i)^*=-q^{\pm N}g^{kh}g_{ki}\partial_h
                                               \label{Real}
\end{equation}
Applying Proposition \ref{propstar} in the latter case
we find that $*$ maps $\c{A}_1$ as well as each of the 
commuting subalgebras $\tilde\c{A}_i^{\pm}$ into itself.

\subsection{Unbraiding `chains' of braided quantum Euclidean 
spaces or spheres}
\label{qspaces}

In this section we consider the braided tensor product of $M\ge 2$ copies
of the quantum Euclidean space $\b{R}_q^N$
~\cite{FadResTak89} (the $\uqs$-covariant quantum space),
i.e.  of the unital associative algebra generated by $x^i$ fulfilling the
relations (\ref{gringo})$_1$,
or of the quotient space of $\b{R}_q^N$ obtained by setting $r^2:=x^ix_i=1$
[the quantum $(N-1)$-dimensional sphere $S_q^{N\!-\!1}$].
(Thus, these will be subalgebras of the Heisenberg algebras
$\c{D}_{+,so(N)},\c{D}_{+,so(N)}$ considered in the
previous subsection).
Again, the multiplet $(x^i)$ carries the fundamental $N$-dimensional 
representation $\rho$ of $\uqs$.
As before, we shall enumerate the different copies of the quantum 
Euclidean space or sphere by attaching an additional Greek index to them, 
e.g.  $\alpha=1,2,...,M$.
The prescription (\ref{absbraiding'}) gives the  cross commutation relations
(\ref{cross})$_1$.

According to Ref. \cite{CerFioMad00}, to define $\varphi_1^{\pm}$ 
(for $q\neq 1$) one actually needs a slightly enlarged version
of $\b{R}_q^N$ (or $S_q^{N\!-\!1}$). One has
to introduce some new generators $\sqrt{r_a}$, with
$0\le a\le \frac N 2$, together with their  
inverses $(\sqrt{r_a})^{-1}$, requiring that 
\begin{equation}
r_a^2=\sum\limits_{h=-a}^a x^hx_h=\sum\limits_{h=-a}^a g_{hk}x^hx^k
\end{equation}
(note that, having set
$n:=\left[\frac N 2\right]$, $r_n^2$ coincides with $r^2$,
whereas for odd $N$ $r_0^2=(x^0)^2$, so we are adding
also $(x^0)^{-1}$ as a new generator). In fact, the commutation relations 
involving these new generators can be fixed consistently, and turn out to
be simply $q$-commutation relations.
$r$ plays the role of `deformed Euclidean distance' 
of the generic `point' $(x^i)$ of $\b{R}_q^N$ from the `origin';
$r_a$ is the `projection' of $r$ on the `subspace' $x^i=0$, $|i|>a$. 
In the previous equation $g_{hk}$ denotes the `metric matrix' of 
$SO_q(N)$,  $g_{ij}=g^{ij}=q^{-\rho_i} \delta_{i,-j}$, which   
is a $SO_q(N)$-isotropic tensor and a deformation of the 
ordinary Euclidean metric. Here,
$(\rho_i):=(n-\frac{1}{2},\ldots,\frac{1}{2},0,-\frac{1}{2},
\ldots,\frac{1}{2}-n)$
for $N$ odd, $(\rho_i):=(n-1,\ldots,0,0,\ldots,1-n)$ for $N$ even. 
$g_{ij}$ is related to the trace projector appearing in (\ref{projectorR})
by $\c{P}_t{}_{kl}^{ij} = (g^{sm}g_{sm})^{-1} g^{ij}g_{kl}$.
The extension of the action of $H$ 
to these extra generators is uniquely determined by the constraints
the latter fulfil. In the case of even $N$ one needs to include also
the FRT generators $\c{L}^+{}^1_1$, $\c{L}^-{}^1_1$ 
(which are generators of $H$) among the generators of $\c{A}_1$. 
In appendix \ref{appeuc} we recall the 
explicit form of $\varphi^{\pm}_1$ in the present case. 
Note that the maps $\varphi^{\pm}_1$ have no analog
in the ``undeformed'' case ($q=1$),
because $\c{A}_1$ is abelian, whereas $H$ is not. 

The unbraiding procedure is recursive. The first step consists of using the
homomorphism $\varphi_1^{\pm}$ found in Ref. \cite{CerFioMad00}
to unbraid the first copy from the others. Following
Theorem \ref{maintheo}, we perform the 
change of generators (\ref{def2.1})$_1$, (\ref{def2.3}) in $\c{A}^-$. 
In view of formula (\ref{imagel-})
we thus find
\begin{equation}
y^{1,i}:=x^{1,i},          \quad\quad\quad
y^{\alpha,i}:=g^{ih}[\mu^1_h,x^{1,k}]_qg_{kj}\,
x^{\alpha,j}, \qquad \alpha>1.                     
                                            \label{def1}
\end{equation}
The suffix 1 in $\mu_a^1$ means that the special elements $\mu_a$
defined in (\ref{defmu}) must be taken as elements
of the first copy. In view of (\ref{defmu})
we see that $g^{ih}[\mu^1_h,x^{1,k}]_qg_{kj}$
are rather simple polynomials in $x^i$ and $r_a^{-1}$,
homogeneous of total degree 1 in the coordinates $x^i$ and $r_a$. 
Using the results given in the appendix
we give now the explicit expression of (\ref{def1})$_2$ for $N=3$:
\begin{eqnarray}
&&y^{\alpha,-} = -qh\gamma_1 \frac r{x^0}x^{\alpha,-} \nn
&&y^{\alpha,0} = \sqrt{q}(q+1) \frac 1{x^0}x^+x^{\alpha,-}+x^{\alpha,0} \\
&&y^{\alpha,+} = \frac{\sqrt{q}(q+1)}{h\gamma_1 rx^0}(x^+)^2 x^{\alpha,-}+
\frac{q^{-1}+1}{h\gamma_1 r}x^+ x^{\alpha,0}- 
\frac 1{qh\gamma_1 r}x^0 x^{\alpha,+} \nonumber
\end{eqnarray}
for any $\alpha=2,...,M$. Here we have set $x^i\equiv x^{1,i}$, 
$h\equiv \sqrt{q}-1/\sqrt{q}$, replaced
for simplicity the values $-1,0,1$ of the indices by the ones $-,0,+$
and denoted  by $\gamma_1\in\b{C}$ a free parameter. 
By Theorem \ref{maintheo}
$y^{1,i}\equiv x^{1,i}$ commutes with $y^{2,i},...,y^{M,i}$.

The alternative unbraiding procedure for the braided tensor
product algebra stemming
from prescription (\ref{absbraiding}) arises by iterating the
change of generators
\begin{equation}
y'{}^{M,i}:=x^{M,i}              \quad\quad\quad
y'{}^{\alpha,i}:=\varphi_M(\c{L}^+{}^i_j)x^{\alpha,j}
=g^{ih}[\bar\mu^M_h,x^{M,k}]_{q^{-1}}g_{kj}\, x^{\alpha,j}, 
                                                   \label{def1'}
\end{equation}
$\alpha<M$.
The special elements $\bar\mu_a$ are defined in (\ref{defmu}), and
suffix $^M$ means that we must take $\bar\mu_a$ as an element
of the $M$-th copy of $\b{R}_q^N$ (or $S_q^{N-1}$).
$y^{M,i}\equiv x^{M,i}$ commutes with  $y^{1,i},...,y^{M-1,i}$.

When $|q|=1$, by a suitable choice (\ref{REAL})
of $\gamma_1,\bar\gamma_1$, as well as of the other
free parameters $\gamma_a,\bar\gamma_a$
appearing in the definitions of $\varphi^{\pm}$
for $N>3$, one can make $\varphi^{\pm}$
into $*$-homomorphisms. 
Applying Proposition \ref{propstar} in the latter case
we find that $*$ maps $\c{A}_1$ as well as each of the 
commuting subalgebras $\tilde\c{A}_i^{\pm}$ into itself.

\appendix
\label{appeuc}

The braid matrix $\hat R$ is related to $\R$ by
$\hat R^{ij}_{hk}\equiv R^{ji}_{hk}:=(\rho^j_h\otimes\rho^i_k)\R$.
With the indices' convention described in sections \ref{qspaces},
\ref{heisenberg} $\hat R$ is given by
\begin{eqnarray}
&&\hat R = q^{-\frac 1N}\left[q \sum_i e^i_i \otimes e^i_i +
\sum_{\scriptstyle i \neq j} e^j_i \otimes e^i_j
+k \sum_{i<j} e^i_i \otimes e^j_j \right]        \label{defRslN} \\   
&&\hat R=q \sum_{i \neq 0} e^i_i \otimes e^i_i +\!
\sum_{\stackrel{\scriptstyle i \neq j,-j} 
{\mbox{ or } i=j=0}}\! e^j_i \otimes e^i_j+ q^{-1} 
\sum_{i \neq 0} e^{-i}_i
\otimes e^i_{-i} \nn
&& \qquad\qquad +k \sum_{i<j}(e^i_i \otimes e^j_j- 
q^{-\rho_i+\rho_j} 
e^{-j}_i \otimes e^j_{-i}) \nonumber
\end{eqnarray}
for $\g=sl(N),so(N)$ respectively.
Here $e^i_j$ is the $N \times N$ matrix with all elements
equal to zero except for a $1$ in the $i$th column and $j$th row.
The braid matrix of $so(N)$ admits the orthogonal projector
decomposition
\begin{equation}
\hat R = q\c{P}_s - q^{-1}\c{P}_a + q^{1-N}\c{P}_t.       \label{projectorR}
\end{equation}
$\c{P}_a,\c{P}_t,\c{P}_s$ are the $q$-deformed antisymmetric, trace,
trace-free symmetric projectors. 
There are just two projectors $\c{P}_a,\c{P}_s$ in
decomposition of the braid matrix of $sl(N)$. the latter is obtained
from (\ref{projectorR}) just by deleting the third term.

We now recall the explicit form of
maps $\varphi^{\pm}$ for the quantum Euclidean spa\-ces
or spheres, found in Ref. \cite{CerFioMad00}.
These are algebra homomorphisms
$\varphi^{\pm}:\b{R}_q^N\cocross U_q^{\pm}so(N)\to \b{R}_q^N$.
We introduce the short-hand notation $[A,B]_x=AB-xBA$. 
The images of $\varphi^-$ (resp. $\varphi^+$)
on the negative (resp. positive) FRT generators
read
\begin{equation}
\varphi^-(\c{L}^-{}^i_j)=g^{ih}[\mu_h,x^k]_qg_{kj}, \qquad\qquad     
\varphi^+(\c{L}^+{}^i_j)=g^{ih}[\bar\mu_h,x^k]_{q^{-1}}g_{kj},
\label{imagel-}
\end{equation}
where
\begin{equation}
\begin{array}{lll}
\mu_0=\gamma_0 (x^0)^{-1}
&\quad\quad\quad\bar\mu_0=\bar\gamma_0 (x^0)^{-1}
&\quad\quad\quad\quad\mbox{for $N$ odd,} \\[6pt]
\mu_{\pm 1}=\gamma_{\pm 1} (x^{\pm 1})^{-1} \c{L}^{\mp}{}^1_1 
&\quad\quad\quad\bar\mu_{\pm 1} = 
\bar\gamma_{\pm 1} (x^{\pm 1})^{-1}\c{L}^{\pm}{}^1_1 
&\quad\quad\quad\quad\mbox{for $N$ even,} \\[6pt]
\mu_a=\gamma_a r_{|a|}^{-1}r_{|a|-1}^{-1} x^{-a}
&\quad\quad\quad\bar\mu_a = 
\bar\gamma_a r_{|a|}^{-1}r_{|a|-1}^{-1} x^{-a}
&\quad\quad\quad\quad\mbox{otherwise,} 
\end{array}                                             \label{defmu}
\end{equation}
and $\gamma_a,\bar\gamma_a \in \b{C}$ are normalization constants fulfilling
the conditions
\begin{equation}
\begin{array}{lll}
\gamma_0 = -q^{-\frac{1}{2}} h^{-1} &\quad\quad
\bar \gamma_0 = q^{\frac{1}{2}} h^{-1}
 &\quad\quad\mbox{for $N$ odd,} \\[6pt]
\gamma_1 \gamma_{-1}=
\left\{\begin{array}{l}
-q^{-1} h^{-2}\\
k^{-2}
\end{array}\right.
&\quad\quad\bar \gamma_1 \bar \gamma_{-1} =
\left\{
\begin{array}{l}
-q h^{-2}\\
k^{-2}
\end{array}\right.
&\quad\quad
\begin{array}{l}
\mbox{for $N$ odd,} \\
\mbox{for $N$ even,}
\end{array}\\[8pt]
\gamma_a \gamma_{-a} =
-q^{-1} k^{-2} \omega_a \omega_{a-1}
&\quad\quad \bar \gamma_a \bar \gamma_{-a}=
-q k^{-2} \omega_a \omega_{a-1}
&\quad\quad\mbox{for $a>1$}. 
\end{array}                                               \label{gamma}
\end{equation}
Here $k:=q\!-\!1/q$, $\omega_a:=(q^{\rho_a}+q^{-\rho_a})$.
Incidentally, for odd $N$
one can choose the free parameters $\gamma_a,\bar\gamma_a$
in such a way that $\varphi^+,\varphi^-$ can be `glued' into
an algebra homomorphism $\varphi:\b{R}_q^N\cocross\uqs\to\b{R}_q^N$
\cite{CerFioMad00}. 
When $|q|=1$, the $*$-structure is given by $(x^i)^*=x^i$
[see (\ref{Real})]. It turns out that
$\varphi^{\pm}$ are $*$-homomorphisms if, in addition, 
\begin{equation}
\begin{array}{ll}
\gamma_{\pm 1}^*=-\gamma_{\pm 1}\qquad & \qquad\mbox{if $N$ even}\\
\gamma_a^*=-\gamma_a\left\{
\begin{array}{l}
1 \mbox{ if }\:\:\:a<0 \\
q^{-2}\mbox{ if }\: a>0
\end{array}\right. 
& \qquad\mbox{otherwise.}
\end{array}
\label{REAL}
\end{equation}

\end{document}